\documentstyle[12pt]{amsart}

\newtheorem{th}{Theorem}[section]	
		\newtheorem{situ}[th]{Situation}
		\newtheorem{prp}[th]{Proposition}

\theoremstyle{definition}       

		\newtheorem{question}[th]{Question}
\newtheorem{pdfn}[th]{Definition}

\theoremstyle{remark}		
\newtheorem{rem}[th]{Remark}

\newcommand{\rar}{\rightarrow} 	
\newcommand{\lrar}{\longrightarrow}	\newcommand{\das}{\dashrightarrow}
\newcommand{\bfp}{{\Bbb{P}}}   	\newcommand{\bfc}{{\Bbb{C}}}
	
\newcommand{\bfz}{{\Bbb{Z}}}	

\newcommand{\co}{{\cal{O}}}	\newcommand{\I}{{\cal{I}}}
\newcommand{\dj}{\cite{dj}}	

\def\cO{{\cal O}}
\def\nor{{\rm nor}}

\def\CHANGE{ }

\newcommand{\spec}{{\operatorname{Spec }}}
\newcommand{\spf}{{\operatorname{Spf }}}

\newcommand{\Aut}{{\operatorname{Aut}}}
\newcommand{\chara}{{\operatorname{char}}}
\newcommand{\Stab}{{\operatorname{Stab}}}

\newcommand{\red}{{\mbox{\small red}}}

\newcounter{cumulative}
\newcommand{\enucum}{\setcounter{enumi}{\value{cumulative}}}
\newcommand{\cumenu}{\setcounter{cumulative}{\value{enumi}}}
\renewcommand{\theenumi}{{\bf\arabic{enumi}}}
\begin{document}

\noindent%
{\large{\bf
Smoothness, Semistability, and Toroidal Geometry\\[2mm]}
\normalsize%
{D. Abramovich\footnote{Partially supported by NSF grant
DMS-9503276.} \\ \small
Department of Mathematics, Boston University\\
111 Cummington, Boston, MA 02215, USA \\
{\tt abrmovic@@math.bu.edu}} \\[2mm]
\normalsize%
{A. J. de Jong
 \\ \small
Department of Mathematics, Harvard University\\
1 Oxford Street, Cambridge, MA 02138, USA \\
{\tt dejong@@math.harvard.edu}}\\[1mm]}
\noindent %
\today
\large
\addtocounter{section}{-1}
\section{INTRODUCTION}
\subsection{Statement} We provide a new proof of the following result:
\begin{th}[Hironaka]\label{resolution}
Let $X$ be a variety of finite type over an algebraically closed field $k$
of characteristic 0, let $Z\subset X$ be a proper closed subset.
There exists a modification $f:X_1 \rar X$,
such that $X_1$ is a quasi-projective nonsingular variety and $Z_1 =
f^{-1}(Z)_\red$
is a strict  divisor of normal crossings.
\end{th}

\begin{rem}
Needless to say, this theorem is a weak version of Hironaka's well known
theorem on resolution of singularities.  Our proof has the feature
that it builds on two standard techniques of algebraic geometry: semistable
reduction for curves,  and toric geometry.
\end{rem}
\begin{rem}
Another proof 
 of the same result was
discovered independently by F. Bogomolov and T. Pantev \cite{bp}.
The two proofs are similar in spirit but quite  different in detail.
\end{rem}

\subsection{Structure of the proof}
\begin{enumerate} \item As in \dj, we choose a projection $X\das P$ of relative
dimension 1, and apply semistable reduction to obtain a model $X'\rar B$ over a
suitable Galois base change $B\rar P$, with Galois group $G$.
\item We apply induction on the dimension to $B$ , therefore we may assume that
$B$ is smooth, and that the discriminant locus of $X'\rar B$ is a $G$-strict
divisor of normal crossings.
\item A few auxiliary blowups make the quotient  $X'/G$ toroidal.
\item Theorem 11* of \cite{te} about toroidal resolutions finishes the
argument.
\end{enumerate}

\subsection{What we do not show}
\subsubsection{Canonicity}
Our proof has the drawback that the resolution
is noncanonical. Some of the steps are not easily carried out in practice,  and
in fact we almost always blow up in the smooth locus. However, it follows
from the proof, that if $C\subset X_{ns}$ is a curve, then we can guarantee
that $X'\rar X$ is an isomorphism in a neighborhood of $C$. We do expect a
slight modification of our argument to give equivariant resolution of
singularities in the case where a finite group acts.

\subsubsection{Positive characteristic} \label{char-p-bad}
There is one crucial point where the proof fails if $\chara\ k =
p>0$. This is at the point where we claim that the quotient $X'/G$ is
toroidal. Given $x\in X'$ and $g\in \Stab\ x$ of order $p$, the action of $g$
on
$\co_{X',x}$ is unipotent. Thus even if $X'$ is toroidal, we cannot guarantee
that the quotient is toroidal as well. It might happen that
the quotient is toroidal by accident, and it would be interesting to see to
what extent such accidents can be encouraged to happen.

In any case, the quotient step goes through if $p \not| \#G$. See remark
\ref{char-p} for a discussion of a bound on $\#G$. As a result, there exists a
function
$$M:\left\{\mbox{\parbox{6cm}{\begin{center}varieties with subvarieties
\end{center}}}\right\} \lrar \bfz$$
 which is bounded on any bounded
family, and which is ``describable'' in a geometrically meaningful way ($M$ for
``multi-genus''), such
that whenever $p>M([X\supset Z])$, our proof goes through for the pair
$X\supset Z$. We were informed by T. Scanlon and E. Hrushovski that the
existence of a function $M$ which is bounded on bounded families is known, by
an application of the compactness theorem in model theory, for any resolution
process, in particular Hironaka's.  A proof of this was given in \cite{ek}.

\subsection{Acknowledgements} We would like to thank B. Hassett, Y. Karshon,
S. Katz, D. Rohrlich, T. Scanlon, and M. Spivakovsky for discussions relevant
to the subject of this paper.

\subsection{Terminology}
We recall some definitions; we restrict ourselves to the case
of varieties over $k$.

A {\bf modification} is a proper birational morphism of irreducible varieties.

An {\bf alteration} $a:B_1\rar B$ is a proper, surjective, generically finite
morphism  of irreducible varieties, see \cite[2.20]{dj}.
The alteration $a$ is a {\bf Galois
alteration} if there is a finite group $G\subset \Aut_B (B_1)$ such that
the associated morphism $B_1/G\rar B$ is birational, compare \cite[5.3]{dj2}.

Let a finite group $G$ act on a (possibly reducible) variety $Z$. Let $Z=\cup
Z_i$ be the decomposition of $Z$ into irreducible components. We say that
{\bf $Z$ is $G$-strict} if the union of translates $\cup_{g\in G} g(Z_i)$ of
each component $Z_i$ is
a normal variety. We simply say that $Z$ is
{\bf strict} if it is $G$-strict for the trivial group, namely every $Z_i$ is
normal.

A divisor $D\subset X$ is called a {\bf divisor of normal crossings} if \'etale
locally at every point it is the zero set of $u_1\cdots u_k$ where
$u_1,\ldots,u_k$ is part of a regular system of parameters. Thus a strict
divisor of normal crossings is what is usually called a divisor of strict
normal crossings,
 i.e., all components of $D$ are nonsigular.

An open embedding $U\hookrightarrow X$ is called a {\bf toroidal embedding} if
locally in the \'etale topology (or classical topology in case $k=\bfc$,
or formally) it is
isomorphic to a torus embedding $T \hookrightarrow V$, (see \cite{te}, II\S
1). If $D=X\setminus U$, we will sometimes denote this toroidal embedding by
$(X,D)$. A finite group action
$G\subset \Aut(U\hookrightarrow X)$ is said to be {\bf toroidal} if the
stabilizer of every point is identified on the appropriate neighborhood with a
subgroup of the torus $T$. We say that a toroidal action is {\bf $G$-strict} if
$X\setminus U$ is $G$-strict. In particular the toroidal embedding itself is
said to be strict if $X\setminus U$ is strict. This is the same as the notion
of {\bf
toroidal embedding without self-intersections} in \cite{te}.

The fundamental theorem about toroidal embeddings we will exploit is the
following:
\begin{th}[\cite{te}, II \S 2, Theorem 11$^*$,
p. 94]\label{mumford}  For any strict
to\-ro\-idal e\-m\-be\-dding
$ U\hookrightarrow X$ 
there exists a canonical
sheaf of ideals $\I  \subset \co_X$ such that the blowup $B_\I(X)$ is
nonsingular.
\end{th}

(See \cite{te}, II \S 2, definition 1 for the notion of canonical
modification.)

\section{The proof}
\subsection{Reduction steps}
We argue by induction on $d=\dim X$. The case $d=1$ is given
by normalizing $X$.

As in \cite{dj}, 3.6-3.10 we may assume that
\begin{enumerate} \item \em $X$ is projective and normal, and
\item $Z$ is the support of a  divisor.
\cumenu
\end{enumerate}

Recall that according to \cite{dj}, lemma 3.11 there is a modification $X'\rar
X$  and a generically smooth
morphism $f:X'\rar P=\bfp^{d-1}$ satisfying the following properties:
\begin{enumerate}
\enucum \it
\item The modification $X'\rar X$  is
isomorphic in a neighborhood of $Z$,
\item every fiber of $f$ is of pure dimension 1,
\item  the smooth locus of $f$ is dense in every fiber, and
\item the morphism $f_{|_{Z}}:Z\rar P$ is finite.
\cumenu
\end{enumerate}
 Since the construction in \dj\ is obtained via a general
projection, we can guarantee that
\begin{enumerate}
\enucum
\item {\it
the generic fiber of $f:X' \rar P$ is
a geometrically connected curve} (this follows from \cite{jou}, 6.3(4), see
\cite{a-pluri}, 4.2).
\cumenu
\end{enumerate}

We now replace $X$ by $X'$.

Applying Lemma 3.13 of \cite{dj}, we may choose a divisor $D\subset X$
mapping
finitely to $P$, which meets every component of every fiber of $f$ in at least
3 smooth points. We may replace $Z$ by $Z\cup D$ (see \cite[3.9]{dj}), and
therefore we have
\begin{enumerate}\it
\enucum
\item $Z$ meets every component of every fiber of $f$ in at least
3 smooth points.
\cumenu
\end{enumerate}

Let  $\Delta_X\subset P$ be the discriminant locus of the map $X\rar P$, and
let
$\Delta_Z\subset P$ be the discriminant locus of the map $Z\rar P$. Let
$\Delta=\Delta_Z\cup \Delta_X$. We can think of $\Delta$ as the discriminant
locus of the pair $Z\subset X\rar P$.

The inductive assumption gives us a resolution of singularities of
$\Delta\subset P$. Thus we may replace $X,Z,P$ and
assume that $P$ is smooth and $\Delta$ is a strict divisor of normal crossings.
Let $\nu:X^\nor\rar X$ be the normalization. By \cite{d-o}, the discriminant
locus
of $\nu^{-1}Z\subset X^\nor\rar P$ is a divisor $\Delta'$ contained in
$\Delta$. We
replace $X$ by $X^\nor$ and $\Delta$ by $\Delta'$. Thus we may assume in
addition to {\bf 1} - {\bf \thecumulative}, that
\begin{enumerate}\it
\enucum
\item $P$ is smooth, and
\item $\Delta$ is a strict normal crossings divisor.
\cumenu
\end{enumerate}

\subsection{Stable reduction}
\renewcommand{\theenumi}{{\bf\alph{enumi}}}

We are now ready to perform stable reduction.
We follow \cite{dj}, 3.18-3.21, but see Remark \ref{opmerking} below.
Let $j:U_0 = P\setminus\Delta\hookrightarrow P$.
Let $a_U: U_0' \rar U_0$ be an \'etale Galois cover
which splits the projection
$Z_{U_0} =j^{-1}(U_0)\cap Z \rar U_0$ into $n$ sections, and trivializes
the 3-torsion subgroup in the
relative Jacobian of $X_{U_0}\rar U_0$. Let $G$ be the
Galois group of this cover. Let $P^\sharp$ be the normalization of $P$ in the
function field of $U_0'$. Let $g$ be the
genus of the generic fiber of $X\rar P$. Let us write
$\overline{{}_3{\bf M}_{g,n}}$ for the Deligne-Mumford compactification
of the moduli scheme of $n$-pointed genus $g$ curves with an
abelian level 3 structure, see \cite{DM},
\cite[2.3.7]{PJ}, \cite[2.24]{dj} and references therein.
There is a `universal'
stable $n$-pointed curve over $\overline{{}_3{\bf M}_{g,n}}$.
We take the closure of the graph of the morphism
$U_0'\rar \overline{{}_3{\bf M}_{g,n}}$ in
$P^\sharp\times \overline{{}_3{\bf M}_{g,n}}$ and obtain a
modification $P'\rar P^\sharp$ and a family of
stable pointed curves $X' \rar P'$. Note that $P'\to P^\sharp$
blows up outside of $U_0'$. We perform a $G$-equivariant blow up
of $P'$ to ensure that we have a morphism $r : X'\to X\times_PP'$,
see \cite[3.18, 3.19 and 7.6]{dj}. Again this blow up has center
outside $U_0'$. In summary:

\begin{situ}\label{situatie}
 There is a  Galois
alteration  $a:P'\rar P$, with Galois group $G$,  and a modification
$r:X'\rar X\times_P P'$ satisfying:
\begin{enumerate} \it
\item the morphism $a$ is finite \'etale over $U_0$,
\item the morphism $r$ is an isomorphism over the open set $U_0\times_P P'$,
\item there are $n$ sections $\sigma_i:P'\rar X'$ such that the proper
transform $Z'$
of $Z$ is the union of their images, and
\item $(X'\rar P',\sigma_1,\ldots,\sigma_n)$ is a stable pointed curve of genus
$g$.
 \cumenu
\end{enumerate}
\end{situ}

\begin{rem}\label{opmerking}
The results of \cite{d-o} imply that the curve $X'$ exists
over the variety $P^\sharp$. We suspect
that the morphism $r$ exists over $P^\sharp$ as well. If this is
true, the following step is redundant.
See \cite{a-pluri} for a similar construction
using Kontsevich's space of stable maps.
\end{rem}

We want to change the situation such that we get $P'=P^\sharp$.
This we can do as follows. Take a blow up $\beta:P_1\to P$ such that
the strict transform $P_1'$ of $P'$ with respect to $\beta$ is finite
flat over $P_1$, see \cite{RG}. By our inductive assumption, we may assume that
$P_1$ is nonsingular and that the inverse image $\Delta_1$ of $\Delta$
in $P_1$ is a divisor with normal crossings. The Galois covering
$U_0'\to U_0$ pulls back to a Galois covering of $P_1\setminus \Delta_1$,
and we get a ramified normal covering $P_1^\sharp$ of $P_1$.
It follows that $P_1^\sharp$ maps to the strict transform $P_1'$.
Hence it is clear that if we replace $P$ by $P_1$, $\Delta$ by $\Delta_1$,
$X'$ by $X'\times_P' P_1^\sharp$, then the morphism $r$ of \ref{situatie}
exists over $P_1^\sharp$.

\renewcommand{\thecumulative}{{\bf\alph{cumulative}}}
We now have in addition to  ({\bf a}) - (\thecumulative),

\begin{enumerate}\it
\enucum
\item $P$ is smooth, $\Delta$ is a strict divisor of normal crossings, and
$a:P'\rar P$ is finite.
\cumenu
\end{enumerate}

Since stable reduction over a normal base is unique (see \cite{d-o}, 2.3)
we have that the action of $G$ lifts to $X'$. Note that the $G$-action
does not preserve the order of the sections.
Note also that since the local fundamental groups of $P\setminus\Delta$
are abelian, so is the stabilizer in $G$ of any point in $P'$. The same
is true for points on $X'$. This property will therefore remain true for
any equivariant modification of $P'$, as well as for points on $X'$.

\subsection{Local description}\label{locdescr}
Denote by $q:P'\rar P$ the quotient map. Let $p\in P $ and let $p'\in
q^{-1}(p)$. Let $s_1,\ldots,s_{d-1}\in \cO_{P,p}$ be a regular system
of parameters on $P$ at $p$ such that $\Delta_x= V(s_1\cdots s_r)$.
Recall that the stabilizer of $p'$ is abelian; this actually
follows from {\bf (a)-(e)} above as the morphism
$P'\rar P$ is ramified only
along the divisor of normal crossings $\Delta$.
Writing $t_i^n = s_i$ for suitable $n$, we
identify a formal neighbourhood of $p'$ in $P'$ as a quotient of the smooth
formal scheme $P''=\spf\ k[[t_1,\ldots t_{r},s_{r+1},\ldots,s_{d-1}]]$ by the
finite
group $(\bfz/n\bfz)^r$ acting by $n$-th roots of unity on the $t_i$.
Thus a formal neighbourhood
of $p'$ in $P'$ is the quotient of $P''$ by a suitable subgroup $H\subset
(\bfz/n\bfz)^r$. \CHANGE This implies that $(\bfz/n\bfz)^r/H$ is identified
with the stabilizer of $p'$ in $G$. \CHANGE

Denote by $X''= P''\times_{P'} X'\rar P''$ the
pulled back stable pointed curve. As
$(\bfz/n\bfz)^r/H \CHANGE\subset G \CHANGE$ acts on $X'$ over
$P'$, we get an action of $(\bfz/n\bfz)^r$ on $X''$ over $P''$.
Let $x\in X''$ be a closed point lying over $p'$, and let $G_x$ be the
stabilizer of $x$ in $(\bfz/n\bfz)^r$. There are two cases:
\renewcommand{\theenumi}{{\bf\roman{enumi}}}
\begin{enumerate}
\item\label{smooth} (Smooth case) Here $x$ is a smooth point of the
morphism $X''\to P''$. In this case the completion of $X''$ at $x$ is
isomorphic to the formal spectrum of the $k[[t_1,\ldots
t_{r},s_{r+1},\ldots,s_{d-1}]]$-algebra
$$k[[t_1,\ldots t_{r},s_{r+1},\ldots,s_{d-1}]][[x]].$$
Here we chose some coordinate $x$ along
the fiber such that $G_x$ acts by a character $\psi_x$ on $x$.
There are two cases with respect to the position of the sections
$Z''\subset X''$:
 \begin{enumerate} \item The point $x\in Z''$. In this case, since $Z''$ is
invariant under the action of $G_x$, we can choose the coordinate $x$ so that
$Z'' = V(x)$.
\item\label{totorify} The point $x\not\in Z''$.  Note that the coordinate $x$
is not
uniquely chosen, and therefore the locus $x=0$ is not uniquely determined.
However, in case $\psi_x$ is nontrivial if we want to see $G_x$ as acting
through the torus for some toroidal structure on $X''$ at $x$
it is necessary to include a locus like $V(x)$ in the boundary.
\end{enumerate}
\item\label{nodal} (Node case) Here $x$ is a node of the fiber of
$X''\to P''$ over $p''$. In this case the completion of $X''$ at $x$ is
isomorphic to the formal spectrum of the $k[[t_1,\ldots
t_{r},s_{r+1},\ldots,s_{d-1}]]$-algebra
$$k[[t_1,\ldots t_{r},s_{r+1},\ldots,s_{d-1}]]
[[x,y]]/(xy - t_1^{k_1}\cdots t_r^{k_r}).$$
We may choose the coordinates $x,y$ such that there
is a subgroup $G_d$ of $G_x$ of index at most $2$ such that $G_d$ acts
by characters on $x,y$, but elements of $G_x$ not in $G_d$
switch the fiber components $x=0$ and $y=0$.
\end{enumerate}

We would like to have the stabilizers acting toroidally on $X''$ in such a way
that the
quotient $X$ becomes strict toroidal.

\subsection{Making the group act toroidally}
There are two issues we need to resolve: in case (\ref{nodal}) above, we want
to modify so that
$G_x/G_d$ disappears from the picture. In case (\ref{totorify}) we
need to modify so that the stratum $x=0$ is not necessary for the toroidal
description of the $G$-action.

In order to describe a global modification we go back to our stable pointed
curve $f:X'\rar P'$.

\subsubsection{Separating branches along nodes}\label{Sepban}
Let $S=\operatorname{Sing} f$ be the singular scheme of the projection $f$.
Let \CHANGE $Y' = B_S(X')\rar X'$ be \CHANGE the blowup of $X'$
along $S$.  Let $Y''$ be the fiber product
$Y'\times_{P'} X''$, and let $y''\in Y''$. We remark that neither
$Y'$ nor $Y''$ is normal in general.

We want to give a local description of $Y'$ and $Y''$. We use the
notation $X_{/x}$ to denote the completion of $X$ at the closed point
$x\in X(k)$, and similar for the other varieties occuring below.
We have $P'_{/p'}=\spf\ R$ where $R$ is the ring
$(k[[t_1,\ldots t_{r},s_{r+1},\ldots,s_{d-1}]])^H$, with $H$ as in
\ref{locdescr} and we
have $X'_{/x'}= \spf\ R[[x,y]]/(xy-h)$ for some monomial $h\in R\subset
k[[t_1,\ldots t_{r},s_{r+1},\ldots,s_{d-1}]]$.
Finally $S=V(x,y)$ scheme theoretically.

{From} this it is easy to read off the following local description,
using that blowing up commutes with completion
in a suitable manner.
\renewcommand{\theenumi}{{\bf\roman{enumi}'}}
\begin{enumerate}
\item\label{smooth'} (Smooth case) $Y''_{/y''}\simeq X''_{/x''}$ as above.
\item\label{nodal'} (Node case) $Y''_{/y''}\simeq \spf\ k[[t_1,\ldots
t_{r},s_{r+1},\ldots,s_{d-1}]][[x,z]]/(xz^2 - t_1^{l_1}\cdots
t_r^{l_r})$. The stabilizer $G_{y''}$ acts diagonally on $t_i, x,z$ (i.e.,
it acts via characters on these elements).
\item\label{double'} (Double case) $Y''_{/y''}\simeq \spf\ k[[t_1,\ldots
t_{r},s_{r+1},\ldots,s_{d-1}]][[x,z]]/(z^2 - t_1^{l_1}\cdots
t_r^{l_r})$. The stabilizer $G_{y''}$ acts
diagonally, but as
in (\ref{totorify}) the coordinate $x$ may not be determined.
\end{enumerate}


The descriptions above determine the local structure of $Y'$ also.
Indeed, $Y'_{/y'}$ is simply the quotient of $Y''_{/y''}$ by the
group $H$ which acts trivially on the coordinate(s) $x$ (and $z$).
We are actually interested in a local description of the normalization
$Y=(Y')^\nor$ of $Y'$.


Let $D\subset Y=(Y')^\nor$ be the union of the
inverse image of $\Delta$ and of $Z$ in $Y$.
{From} the local descriptions given above it is already clear
that $Y\setminus D \hookrightarrow Y$ is a strict toroidal embedding.
Moreover, $D$ is $G$-strict, since $\Delta$ is
strict and since the blowup $Y\rar X'$ separates fiber components.
However, the action of $G$ on the pair $(Y\setminus D, Y)$
is not yet toroidal: indeed, in case (\ref{smooth'}) and
in case (\ref{double'}) 
if the character $\psi_x$ is nontrivial, we have a problem.
More precisely, this is the situation explained in (\ref{totorify}).

\subsubsection{Torifying a pre-toroidal action}\label{Torlta}
Here we show how to do one canonical blow up $Y_1\rar Y$
(analogousely to \cite{te}, II \S 2) which makes the action of $G$ toroidal.
The situation \CHANGE $(Y, D, G)$ we reached at the end of
Section \ref{Sepban} is summarized by the conditions in Defintion
\ref{pretoroidal} below. \CHANGE

\begin{rem}
We expect that this discussion should be of interest in a more general
setting.
\end{rem}

\begin{pdfn}\label{pretoroidal} Let $U=Y\setminus D\subset Y$ be a toroidal
embedding, $G\subset
\Aut(U\subset
Y)$ a finite subgroup, such that $D$ is $G$-strict. For any point $y\in Y$
denote the stabilizer of $y$ by $G_y$. We say that the action of $G$ is {\bf
pre-toroidal} if at every point $y\in Y$, either $G_y$ acts toroidally at $y$,
or the following situation holds:
\begin{itemize}
\item There exists an isomorphism $\epsilon: Y_{/y}\cong \spf\ R[[x]]$,
\item where $\spf R$ is the completion of a toroidal embedding
$T_0\hookrightarrow Y_0$ at a point $y_0\in Y_0$,
\item where $D_{/y}$ corresponds to
$(Y_0\setminus T_0)_{/y_0}\times_{\spf\ k} \spf\ k[[x]]$,
\item where $G_y$ acts toroidally 
on $T_0\hookrightarrow Y_0$ fixing $y_0$
and
\item $G_y$ acts on the coordinate $x$ via a character $\psi_x$ such that the
isomorphism $\epsilon$ of completions is $G_y$-equivariant. In other words,
$G_y$ acts
toroidally on  $$(T_0)_{/y_0}\times_{\spf\ k}\spf\ k[[x, x^{-1}]]\quad \subset
\quad
 Y_{/y}.$$
\end{itemize}
\end{pdfn}

Analogousely to definition 1 of \cite{te}, II \S 2, we define pre-canonical
ideals:

\begin{pdfn}
Let $G\subset \Aut(U\subset Y)$ be a pre-toroidal action. A $G$-equi\-va\-riant
ideal sheaf
$\I\subset \co_Y$ is said to be {\bf pre-canonical} if the following holds:

 For any $y, y'$ lying on
the same stratum, and any isomorphism $\alpha: O_{Y_{/y}} \rar O_{Y_{/y'}}$
preserving the strata and inducing an isomorphism carrying $G_y$ to $G_{y'}$,
we have $\alpha(\I_{/y}) =\I_{/y'}$.

If $\I$ is pre-canonical, we say that the normalized  blowup
$b:(B_\I Y)^\nor\rar Y$ is a pre-canonical blowup.
\end{pdfn}

\begin{pdfn}
We say that a pre-canonical blowup $\tilde{Y} \rar Y$ {\bf torifies} $Y$ if
$G$ acts toroidally on $(b^{-1} U \subset \tilde{Y})$.
\end{pdfn}

\begin{th}\label{opblazen}Let $G\subset \Aut(U\subset Y)$ be a pre-toroidal
action. Then there
exists a canonical choice of a pre-canonical ideal sheaf $\I_G$ such that the
pre-canonical blowup $b:(B_{\I_G} Y)^\nor\rar Y$ torifies $Y$.
\end{th}

The theorem follows immediately from the affine case below:

\begin{prp}\label{affineopblazen}
Let $T_0\subset X_0$ be an affine torus embedding, $X_0 = \spec\
R$.
Let $G\subset T_0$ be a finite subgroup of $T_0$, let $p_0\in X_0$ be
a fixed point of the action of $G$, and let $\psi_x$ be a character
of $G$. Consider the torus embedding of $T=T_0\times \spec\ k[x,x^{-1}]$
into $X= X_0\times \spec\ k[x]$, where we let $G$ act
on $x$ via the character $\psi_x$. Assume that the map
$G\to T$ induced from this is injective. Write $p=(p_0,0)\in X$ and
write $D=(X_0\setminus T_0) \times \spec\ k[x]$.
There is a canonical $T$-equivariant ideal $I_G\subset R[x]$,
satisfying the following:

Let $b:X_1 = (B_{I_G}X)^\nor\rar X$, the normalization of the blowup of $X$
along $I_G$. Let $U_1=b^{-1}(T_0\times \spec\ k[x])$. Then
$U_1\hookrightarrow X_1$ is toroidal and $G$ acts toroidally on
$U_1\hookrightarrow X_1$.

If $X_0', T_0', G', p_0'$ and $\psi'_x$ is a second set of such data,
and if we have an isomorphism of completions
$$ \varphi : X_{/p}\cong X'_{/p'}, $$
which induces isomorphisms $G\cong G'$ and
$D_{/p}\cong D'_{/p'}$, then
$\varphi$ pulls back $I_G$ to the ideal $I_{G'}$.

Furthermore, if $q_0$ is any point of $X_0$ and if $G_q\subset G$ is
the stabilizer of $q$ in $G$, then the stalk of $I_G$ at $q$
is the same as the stalk of $I_{G_q}$ at $q$.

\end{prp}

\begin{rem} The ideal $I_G$ is
called the {\bf torific ideal} of the situation $G\subset  \Aut(T_0\times
\spec\ k[x] \subset X)$.
\end{rem}

{\bf Proof.}
 For any monomial  $t\in R[x]$ let $\chi_t $ be the associated
character. We restrict these characters to $G$ and obtain a character
$\psi_t: G\rar k^\ast$.

Define $M_x=\{t | \psi_t = \psi_x\}$, the set of mononials on which
$G$ acts as it acts on $x$. Notice that the character $\psi_x$ of $G$
is uniquely determined by the data $G\to \Aut(X_{/p},D_{/p})$.
Define the torific ideal $I_G= \left< M_x\right >$, the ideal generated by
$M_x$. We see that it satisfies the variance property of the
proposition with respect to isomorphisms $\varphi$. We leave
the localization property at the end of the proposition as an
excercise to the reader.

Define $b:X_1 \rar X$ as in the proposition.

%
%

Since $G$ is a subgroup of $T_0$, we have that $M_x$ contains a monomial in
$R$. Therefore the ideal $I_G$
is generated by $x$ and a number of monomials $t_1,\ldots,t_m\in R\cap M_x$.
The blow up has a chart associated to each of the generators
$x, t_1,\ldots, t_m$. On the chart ``$x\neq 0$'' we have that  the
inverse image of $D$ contains the inverse image
of $D\cup V(x)$, and hence the action of $G$ is toroidal, being
toroidal with respect to $T\subset X$ on $X$. The other charts ``$t_i\neq 0$''
can be described as the spectra of the rings
$$\tilde R= R[x][u, \{s_j\}_{j\not=i}]/(ut_i-x, s_jt_i-t_j, h_\alpha),$$
that is $u = x/t_i$ and $s_j =t_j/t_i$. Since there are no
relations between $x$ and $t_j$ we can take the $h_\alpha$ to be certain
polynomials in $R[s_j]$  (the ideal
generated by the $t_i$ is flat over $k[x]$). Note that $G$ fixes the element
$u$ in this algebra.

Thus it follows that the normalization of the ring $\tilde R$
is of the form $R'[u]$, where $G$ acts trivially on $u$, $R'$ is
the ring associated to an affine torus embedding and $G$ acts toroidally
on $\spec\ R'$.
\qed

\subsubsection{Strictness}\label{striktheid}
We return to the triple $(Y, D, G)$ we obtained at the end of \ref{Sepban},
in particular we have the $G$-equivariant map $Y\to P'$ and the Galois
alteration $P'\to P$ with group $G$.
Denote by $b: Y_1 \rar Y$ the torifying blowup obtained by normalizing the
blowup at the torific ideal of $Y$, as in Theorem \ref{opblazen}.
It remains to check that the
divisor $b^{-1}(D)\subset Y_1$ is $G$-strict, so that the quotient is a
strict toroidal embedding.

First, although a coordinate $x$ for the pre-toroidal action is not globally
defined, we may always find such coordinates on Zariski neighborhoods of the
relevant points. Let $y\in Y$ be such a point. Let
$\epsilon: Y_{/y}\cong \spf\ R[[x]]$ be as in Defintion \ref{pretoroidal}.
Let $x'\in O_{Y,y}$ be an element of the local ring of $Y$ at
$y$ that is congruent to $\epsilon^{-1}(x)$ up to a high power of the
maximal ideal. Then the element
$$x''=\sum_{g\in G_y} g(x') \psi_x^{-1}(g) / |G_y| \in O_{Y,y}$$
transforms according to $\psi_x$ under the action of $G_y$ and
is congruent to $\epsilon^{-1}(x)$ up to a high power of the
maximal ideal. We may then change the isomorphism $\epsilon$ so
that the element $x''$ will correspond to the coordinate $x$.
Therefore, we may assume that there is a $G_y$-invariant
Zariski open neighbourhood $W=W_y$ of $y$ and a function
$x\in \Gamma(W, \cO)$ that transforms according to the character
$\psi_x$ under the group $G_y$, and giving rise to the local
coordinate for a suitably chosen isomorphism $\epsilon$ as
in Definition \ref{pretoroidal}. After possibly shrinking
$W$, we get that $W$ is strict toroidal with respect
to the divisor $D_W= (W\cap D) \cup V(x)$. If we choose  $W$ sufficiently
small, we may assume that the conical polyhedral complex of $(W, D_W)$ a
single cone $\sigma_W$ (see \cite{te}, II \S 1 definition 5, through p. 71).

Now, to show that  $b^{-1}({D})$ is strict, it suffices to show this on
an affine neighborhood of any point. This follows immediately from
theorem $1^*$ of \cite{te}, applied to $(W,D_W)$.

At this point it should be remarked that we could get away without
$G$-strictness: using the constructions in \cite{te}, I \S2 lemmas 1-3 on pages
33-35, and \cite{dj}, 7.2, it is not difficult to construct a $G$-equivariant
blowup which is a $G$-strict toroidal embedding. Still it is of interest to
know that $b^{-1}(D)\subset Y_1$ is already $G$-strict.

 Let $E\subset b^{-1}(D)$ be a component
and let $g\in G$ be an element such that $g(E)\cap E\not= \emptyset$.
We have to show that $g(E)=E$. In the case that $b(E)$ is a component
of $D$, this follows immediately from the $G$-strictness of $D$,
see end of \ref{Sepban}; the difficult case is that of a component
that is collapsed under $b$.

Let $y$ be a general point of the intersection $g(E)\cap E$,
and let $y=b(y)$. Let $D_1,\ldots, D_s$ be the components of
$D$ at $y$.
The components $g^{-1}(D_i)$ are the components of
$D$ at $g^{-1}(y)$.
After reordering, we may assume that $D_1,\ldots, D_r$ are the components
which contain $b(E)$, for some $r\leq s$. As $y\in g(E)\cap E$ we
get $y\in g(b(E))\cap b(E)$, hence $y\in
g(D_i)\cap D_i$
for all $i=1,\ldots, r$. By $G$-strictness of $D$, we get
$g(D_i)=D_i$ for $i=1,\ldots,r$.

Let $\sigma_y$ be the cone corresponding to the
toroidal structure $D_{W_y}$ at $y$ of the second paragraph above.
The divisors $D_i$ correspond to rays $\tau_{D_i,y}$ in the boundary of
$\sigma_y$. The divisor $E$ corresponds to a
ray $\tau_{E,y}$ in $\sigma_y$, belonging to the polyhedral
decomposition associated to the given blowup.
For the Zariski open neighbourhood $W_{g^{-1}(y)}$ we take
$g^{-1}(W_y)$ and for the divisor $D_{W_{g^{-1}(y)}}$ we take
$g^{-1}(D_{W_y})$. Hence we get an identification
$\sigma_y\cong \sigma_{g^{-1}(y)}$ given by $g^{-1}$.
Under this identification the ray $\tau_{E,y}$ is mapped
to the ray $\tau_{g^{-1}(E), g^{-1}(y)}$. Note that
$E$ passes through $g^{-1}(y)$, hence we have the
ray $\tau_{E, g^{-1}(y)}$ in $\sigma_{g^{-1}(y)}$.
We have to show that $\tau_{E,g^{-1}(y)}=\tau_{g^{-1}(E), g^{-1}(y)}$.

Let $M_y$ be the group of Cartier divisors supported along $D_{W_y}$ (see
\cite{te}, II \S 1, definition 3). After
shrinking $W_y$ we may assume that $M_y$ is the free abelian group generated by
the divisors of
rational functions
$f_j\in \Gamma(W_y,\cO), j=1,\ldots,m$ and the function $x$ on $W_y$.
 We may consider
$f_j$ and $x$ as rational functions on $Y_1$ as well.
The ray $\tau_{E,y}$ is determined by the order of vanishing
of the functions $f_j$ and $x$ along $E$ on $Y_1$,
let us call these orders $n_j$ and $n$.
Similarly, $\tau_{E,g^{-1}(y)}$ (resp.\ $\tau_{g^{-1}(E), g^{-1}(y)}$)
is determined by the order of vanishing of the functions
$g^*(f_j)$ and $g^*(x)$ along $E$ (resp.\ $g^{-1}(E)$), let
us call these orders $n'_j$ and $n'$ (resp.\
$n''_j$ and $n''$).

It is clear that $n_j=n''_j$ and $n=n''$.
Notice that $g^*(f_j)$ are supported along $g^{-1}(D_i)$ in $W_{g^{-1}(y)}$.
Here $g^{-1}(D_i) = D_i$ for $i\leq r$, so that  $f_j/g^*f_j$
has order 0 along $D_i, i\leq r$.
Since the vanishing along $E$ may be computed on
$b^{-1}(W_y\cap W_{g^{-1}(y)})$, we see that $n_j=n_j'$.
We are through if we show that $n=n'$; this is equivalent to
showing that the order of vanishing of $x$ along
$E$ is the same as the order of vanishing of $g^{-1}(x)$ along $E$.

The component $E$ is exceptional for the morphism $b$.
The torifying blowup $b$ over $W_y$ blows up inside $V(x)$
only. Hence we see that $b(E)\cap W_y$ is contained in
$V(x)$, i.e., $n>0$. The same argument applied to $g^{-1}(x)$
on $W_{g^{-1}(y)}$ works to show that $n'>0$. For a general
point $z\in b(E)$, we see that both $x$ and $g^{-1}(x)$
define a local coordinate that can be used to define the
local toroidal structure around $z$, as in Definition
\ref{pretoroidal}. Therefore, we need only to show that
the valuation of $x$ along $E$ does not depend on the choice
of the parameter $x$ as in Definition \ref{pretoroidal}.

Let $\epsilon: Y_{/z}\cong \spf\ R[[x]]$ be as in
Definition \ref{pretoroidal}. Any other $x'\in R[[x]]$
that gives a local coordinate for some
pretoroidal structure is of the form
$$x'= u x + r_0 + \sum_{j\geq 2} r_j x^j, $$
with $u, r_0, r_j\in R$ and where $u$ is a unit and
$r_0$ is in the maximal ideal of $R$. The element
$u$ being a unit,
we may divide by it without changing the
orders of vanishing on $Y_1$. Thus we may assume $u=1$.
Consider the family of automorphisms $\varphi_t$,
$t\in [0,1]$ of $R[[x]]$ given by $\varphi_t(r)=r$, $r\in R$ and
$$ \varphi_t(x)=x+t\Big(r_0+ \sum_{j\geq 2} r_j x^j\Big).$$
These are automorphisms that occur in Proposition
\ref{affineopblazen}. Hence these act on the
formal completion $Y_1^\wedge$ of $Y_1$ along $b^{-1}(z)$,
in view of Proposition \ref{affineopblazen}.
Since they form a continuous family (acting continuously
on the charts described in the proof of Proposition
\ref{affineopblazen}), they will fix (formal) components
of $b^{-1}(D_{/z})$ such as
$E^\wedge$ in $Y_1^\wedge$. Hence the orders
of vanishing of the functions $\varphi_t(x)$ along
the (formal) component $E^\wedge$ are all the same.
In particular, we get the equality for $x$ and $x'$, as
desired.

\subsection{Conclusion of proof}
If we combine the results of \ref{Sepban}, \ref{Torlta}
and \ref{striktheid} then we see that we may assume our
Galois alteration $(Y_1, D_1, G)$
of $(X,Z)$ is such that $Y\setminus D_1\hookrightarrow Y_1$ is a $G$-strict
toroidal embedding.
Hence the quotient $(Y_1\setminus D_1)/G\hookrightarrow Y_1/G$ is a
strict toroidal embedding, and $Y_1/G\to X$ is a modification.
Hence we may replace $(X, Z)$ by $(Y_1/G, D_1/G)$.
By Theorem \ref{mumford}, there exists a (toroidal)
resolution of singularities of this pair.
This ends the proof of Theorem \ref{resolution}.
\qed

\begin{rem}\label{char-p}
If $\chara\ k = p>0$ the proof goes through if it works on $P$, and if the
exponent of the Galois group $G$ is small enough. Since $G$ can be taken as the
Galois group of the torsion on a generalized Jacobian, we can bound it in terms
of  $g+n$, where $g$ is the relative genus and $n$ is the degree of the divisor
$Z$. Given a family of varieties of finite type $X\rar S$, we
can
make the constructions in the proof uniform over $S$. Roughly speaking,
after replacing
$S$ by a dense open, there is a sequence of projections $X\das
X_1\das\cdots\das S$ and relative divisors $D_i \subset X_i$ which will do the
job. Therefore  the order of all groups involved is bounded in terms of the
relative genus of $X_i\das X_{i+1}$ and the degree of $D_i\das X_{i+1}$. Thus
there is a ``geometrically meaningful'' function $M$, as
described in \ref{char-p-bad}.
\end{rem}


\begin{thebibliography}{HHHHHHH}

\bibitem[$\aleph$]{a-pluri} D. Abramovich,
{\em Galois pluri-nodal reduction,} \S 4 in preprint.

\bibitem[B-P]{bp} F. Bogomolov and T. Pantev, {\em Weak Hironaka theorem},
preprint.

\bibitem[DM]{DM} P. deligne and D. Mumford,
{\em The irreducibility of the space of curves of given genus},
Publications Math\'ematiques I.H.E.S. {\bf 36}, 75-109 (1969).

\bibitem[E]{ek} P. Eklof, {\em Resolutions of singularities in prime
characteristics for almost all primes}, Trans. AMS {\bf 146} (1969),
429-438.

\bibitem[H]{h}  H. Hironaka, {\it Resolution of singularities of an
algebraic variety over a field of characteristic zero: I, II,\/} Ann. of Math.
(2) {\bf 79} (1964), 109-326.

\bibitem[dJ]{dj} A. J. de Jong, {\em Smoothness, semistability, and
alterations}, preprint.\\
{\tt ftp://ftp.math.harvard.edu/pub/AJdeJong/alterations.dvi}

\bibitem[dJ2]{dj2} A. J. de Jong, {\em Families of curves and
alterations}, Preprint.\\
{\tt ftp://ftp.math.harvard.edu/pub/AJdeJong/curves.dvi}

\bibitem[dJ-O]{d-o} A. J. de Jong and F. Oort, {\em On extending families of
curves}, preprint.\\
{\tt ftp://ftp.math.harvard.edu/pub/AJdeJong/extend.dvi}

\bibitem[J]{jou} J.-P. Jouanolou, {\em Th\'eor\`emes de Bertini et
Applications}, Progress in Math. 42, Birkh\"auser: Boston, Basel, Stuttgart,
1983.

\bibitem[KKMS]{te} G. Kempf, F. Knudsen, D. Mumford and B. Saint-Donat,
{\em Toroidal Embeddings I}, Springer, LNM 339, 1973.

\bibitem[PJ]{PJ} M. Pikaart and A. J. de Jong,
{\em Moduli of curves with non-abelian level structures,}
in: The moduli spaces of curves, R. Dijkgraaf, C. Faber,
G. van der Geer (ed.), Birkh\"auser, Boston-Basel-Berlin,
Progress in math. vol. 129, 483-509 (1995).

\bibitem[RG]{RG} M. Raynaud and L. Gruson,
{\em Crit\`eres de platitude et de projectivit\'e,}
{Techniques de $\ll$platification$\gg$ d'un module},
Inventiones mathematicae {\bf 13}, 1-89 (1971).


\end{thebibliography}
\end{document}